\documentclass{article}%
\usepackage{amssymb}
\usepackage{amsfonts}
\usepackage{amsmath}
\usepackage{graphicx}%
\providecommand{\U}[1]{\protect\rule{.1in}{.1in}}
\newtheorem{theorem}{Theorem}

\newtheorem{corollary}[theorem]{Corollary}

\newtheorem{lemma}[theorem]{Lemma}

\begin{document}

\title{A non-iterative method for electrical property\ tomography based on a simple formula}
\author{Victor P. Palamodov\\Tel Aviv University}
\date{}
\maketitle

\textbf{Abstract }A non-iterative method of reconstruction is proposed from
data of MRI system and of a harmonic electro-magnetic field at Larmor
frequency. The method is based on the exact analytic formula for the contrast
source function. A geometric method for acquisition of the full inductive
field is discussed.



\section{Introduction}

Electrical properties tomography (EPT)\ is a noninvasive reconstruction
technique for retrieving electrical properties (conductivity and permittivity)
of biological tissue from magnetic fields generated by radio frequency coils
in a magnetic resonance imaging scanner. The electric properties of tissue are
of great interest, since these properties can be used to aid the
discrimination of cancerous tissue from benign tissue and characterize various
kinds of pathological tissue of human body \cite{3}, \cite{5}. The main
benefits of EPT over other reconstruction modalities is that it uses the
Larmor frequency fields of an MRI system, which can penetrate biological
tissues and it does not make use of surface electrode mounting, current
injection or additional hardware.

The method called now EPT was proposed by Haacke \textit{et al}. \cite{1}. The
first application of this method was done by Wen \cite{2} who used an
approximation for the contrast source function $\gamma/\gamma_{0}$ (CSF). Song
and Seo \cite{10} reduced reconstruction of\ admittivity $\gamma$\ to solution
of an elliptic equation under assumption $\partial\gamma/\partial z=0$ in
terms of $\mathit{B}^{+}$ field. See also \cite{6} for an explicit
construction. Ammari \textit{et al \cite{12} }studied the problem in terms of
the quasilinear system of differential equations. In \cite{7} a survey of
local methods of reconstructions is given.\ C van den Berg \textit{et al}
\cite{8} developed the iterative method of determination of the CSF\ of an
object making use of a global integral approach. Arduino \textit{et al
\cite{Ard} }applied the iterative conjugate gradient method (1000 steps). We
describe here a non-iterative\ algorithm\textit{ }for retrieving CSF which is
based on the exact global formula. This is the first algorithm of this kind so far.

\section{Basic equations}

We assume that the conductivity and permittivity are isotropic at the angular
Larmor frequency $\omega$ and ignore spatial permeability variations $\mu
=\mu_{0}$, since they are considered small for biological tissue. Let $\Omega$
be an object domain in the MRI scanner that is disjoint of the antenna (coil)
generating the RF wave.$\ $The total time harmonic electromagnetic field in
$\mathbb{R}^{3}$ can be written in the form $\left\{  E\exp\left(  i\omega
t\right)  ,H\exp\left(  i\omega t\right)  \right\}  \ $where%
\[
E=\mathbf{E}+E^{sc},\ H=\mathbf{H}+H^{sc},
\]
$\left\{  \mathbf{E=E}^{inc},\mathbf{H=H}^{inc}\right\}  $\ is the incident
and $\left\{  E^{sc},H^{sc}\right\}  \ $the\ scattered\ fields$\ $due to the
presence of the object. Following to van der Berg \textit{et al }\cite{5} we
use the integral representations%
\begin{equation}
E^{sc}=\left(  \mathbf{k}^{2}+\nabla\nabla\cdot\right)  A, \label{2}%
\end{equation}%
\begin{equation}
H^{sc}=\gamma_{0}\nabla\times A \label{3}%
\end{equation}
for the scattered fields where $\gamma_{0}=i\omega\varepsilon_{0}%
$,$\ \varepsilon_{0}$ is the permittivity in vacuum and%
\begin{equation}
A\left(  p\right)  =\int_{\Omega}G\left(  p-q\right)  \chi\left(  q\right)
E\left(  q\right)  \mathrm{d}V\left(  q\right)  ,\ q,p\in\mathbb{R}^{3}
\label{4}%
\end{equation}
is the vector potential and
\begin{equation}
G\left(  q\right)  =-\frac{\exp\left(  ik\left\vert q\right\vert \right)
}{4\pi\left\vert q\right\vert },\ k=\frac{\omega}{c_{0}} \label{14}%
\end{equation}
is the scalar Green's function for\ the Helmholtz operator $\Delta+k^{2}$ in
3D\ background medium. The contrast source function is defined as%
\[
\chi\left(  q\right)  =\frac{\gamma}{\gamma_{0}}-1
\]
where admittivity $\gamma=\sigma+i\varepsilon\omega$ is the per meter
admittance. Conductivity $\sigma$ and permittivity $\varepsilon$ are profiles
of the object, $\Omega\ $\ is a domain that contains the support of $\chi.$

\section{Evaluation of the contrast function}

We write the electric and magnetic\ fields as first order differential forms
\[
E=E_{x}\mathrm{d}x+...,\ B=B_{x}\mathrm{d}x+...
\]
Hodge star operator $\ast$ acts linearly on differential forms $\alpha,\beta$
on\textrm{ }$\mathbb{R}^{3}$\textrm{ }by the rule (\cite{9}, p.15)%
\begin{align*}
\alpha^{\ast}\wedge\beta &  =\alpha\wedge\beta^{\ast}=\left\langle
\alpha,\beta\right\rangle \mathrm{d}x\wedge\mathrm{d}y\wedge\mathrm{d}z,\\
\left(  a\mathrm{d}x\right)  ^{\ast}  &  =a\mathrm{d}y\wedge\mathrm{d}%
z,\ \left(  a\mathrm{d}y\wedge\mathrm{d}z\right)  ^{\ast}=a\mathrm{d}%
x;\ \left(  x\rightarrow y\rightarrow z\rightarrow x\right)  ,
\end{align*}
where $\alpha,\beta\ $are arbitrary differential forms in $\mathbb{R}^{3}%
\ $of\ equal$\ $degree and $\left\langle ,\right\rangle \ $means the natural
scalar product of the forms. The dual differential is defined by
$\mathrm{d}^{\ast}=\ast\mathrm{d}\ast$ and Laplace operator $\Delta
=\mathrm{d}^{\ast}\mathrm{d}+\mathrm{dd}^{\ast}$ commutes with $\mathrm{d}%
=\nabla$ and $\mathrm{d}^{\ast}.\ $

\begin{theorem}
Let $\Omega$\ be a domain in $\mathbb{R}^{3}$ disjoint of the RF\ antenna. If
the fields $\mathbf{E}$ and $B$\ are known then the contrast source function
can be found in $\Omega$ by%
\begin{equation}
\chi=-\frac{\left(  \mathbf{E}-\ast\mathrm{d}B\right)  \wedge\left(
\Delta+k^{2}\right)  B^{\ast}}{\left(  \mathbf{E}-\ast\mathrm{d}B\right)
\wedge\left(  \mathrm{d}\mathbf{E}+k^{2}B^{\ast}\right)  } \label{9}%
\end{equation}
if the denominator does not vanish.
\end{theorem}

\begin{corollary}
We have%
\begin{equation}
\frac{\gamma}{\gamma_{0}}=\chi+1=\frac{\left\langle \mathbf{-}\Delta B^{\ast
}+\mathrm{d}\mathbf{E},\ast\mathrm{d}B-\mathbf{E}\right\rangle }{\left\langle
k^{2}B^{\ast}+\mathrm{d}\mathbf{E},\ast\mathrm{d}B-\mathbf{E}\right\rangle }.
\label{15}%
\end{equation}

\end{corollary}

Note that $\left\langle \mathrm{d}\mathbf{E,E}\right\rangle =0$ since
$\mathrm{d}\mathbf{E}=-i\omega c^{-1}\mu\mathbf{H}^{\ast}$\ and $\mathbf{H}%
^{\ast}\wedge\mathbf{E}=\left\langle \mathbf{H,E}\right\rangle =0.\ $

\textbf{Remark. }Equation
\begin{equation}
\gamma=\frac{1}{i\omega\mu_{0}}\frac{\left\langle \nabla^{2}H,\nabla\times
H\right\rangle }{\left\langle H,\nabla\times H\right\rangle } \label{12}%
\end{equation}
is mentioned\textit{ }by Seo\ \cite{10} and attributed to Nachman \textit{et
al} \cite{11}.

\textit{Proof of Theorem 1}\textbf{. }

\begin{lemma}
\label{P}Let $\Omega$ be a compact set in $\mathbb{R}^{n}$ that can be
contracted to a point in itself.\ For any $k\geq1,\ $and arbitrary closed
differential $k$-form $a$ on $\Omega$, there exists a $k-1$-form $b$ on
$\Omega\ $such\ that$\ \mathrm{d}b=a\ $on $\Omega.$
\end{lemma}

For a proof see "Converse of the Poincar\'{e} lemma" \cite{9},
p.29.\footnote{We may assume that coefficients of $a$ and of $b\ $are
distributions since $b$ will not appear in the final formula.}.

Note that equation (\ref{9}) does not depend on $\Omega.\ $Therefore we may
assume that$\ \Omega$ is a compact set with smooth boundary in $\mathbb{R}%
^{3}\ $that fulfils the condition of Lemma \ref{P}.

\begin{lemma}
\ The system of equations for a 1-form $C$%
\begin{equation}
\mathrm{d}C=B^{\ast},\ \mathrm{d}^{\ast}C=0 \label{11}%
\end{equation}
has a solution defined on $\Omega.\ $There exists a function $g$ on $\Omega$
such that
\begin{equation}
C+\mathrm{d}g=A_{0}. \label{13}%
\end{equation}

\end{lemma}

\textit{Proof.} Form $B^{\ast}$ has bounded coefficients and fulfils the Gauss
law $\mathrm{d}B^{\ast}=0$ on $\Omega.\ $By Lemma \ref{P} there exists 1-form
$C_{0}$ on $\Omega\ $that satisfies $\mathrm{d}C_{0}=B^{\mathbf{\ast}}$.
Consider Dirichlet problem%
\[
\Delta f=-\mathrm{d}^{\ast}C_{0}%
\]
for a function $f$ on $\Omega.\ $To solve it we define the function $h$ on
$\mathbb{R}^{3}$ that is equal to $-\mathrm{d}^{\ast}C_{0}\ $on $\Omega$ and
$h=0$ on the complement to $\Omega.$\ Set $f=G_{0}\ast h$\ where $G_{0}$ is
the kernel (\ref{14}) with $k=0.\ $Differential form $C=C_{0}+\mathrm{d}f$
fulfils (\ref{11}) since%
\[
\mathrm{d}^{\ast}C=\mathrm{d}^{\ast}C_{0}+\mathrm{d}^{\ast}\mathrm{d}%
f=\mathrm{d}^{\ast}C_{0}+\Delta f=0.
\]
Equation (\ref{3}) implies%
\[
B^{\ast}=\mu\left(  H^{sc}\right)  ^{\ast}=\frac{i\omega}{c^{2}}%
\mathrm{d}A=\mathrm{d}A_{0},
\]
where $A_{0}=i\omega c^{-2}A$ and $A$ is the vector potential (\ref{4}).\ By
(\ref{11}) we have $\mathrm{d}\left(  A_{0}-C\right)  =0.$ Again, by Lemma
\ref{P} there exists a function $g$ on $\Omega$ such that $A_{0}%
-C=\mathrm{d}g$ and (\ref{13}) follows.$\ \blacktriangleright$

By (\ref{4}) and (\ref{2}) we have%
\[
A=G_{\chi}E=G_{\chi}\left(  E^{sc}+\mathbf{E}\right)  =G_{\chi}\left(
k^{2}+\mathrm{dd}^{\ast}\right)  A+G_{\chi}\mathbf{E}%
\]
since $\nabla\nabla\cdot\ =\mathrm{dd}^{\ast}\ $where$\ G_{\chi}u=G\ast\left(
\chi u\right)  $. This yields%
\[
G_{\chi}\mathbf{E}=A-G_{\chi}E^{sc}=\left(  I-G_{\chi}\left(  k^{2}%
+\mathrm{dd}^{\ast}\right)  \right)  A_{0}.
\]
By (\ref{13}) we have
\[
G_{\chi}\mathbf{E}=\left(  I-G_{\chi}k^{2}\right)  C+\left(  I-G_{\chi}\left(
k^{2}+\Delta\right)  \right)  \mathrm{d}g
\]
since $\nabla\nabla\cdot\mathrm{d}g=\mathrm{dd}^{\ast}\mathrm{d}%
g=\Delta\mathrm{d}g.\ $Apply Helmholtz operator to both parts and get%
\begin{align*}
\chi\mathbf{E}  &  =\left(  k^{2}+\Delta-\chi k^{2}\right)  C+\left(  \left(
k^{2}+\Delta\right)  -\chi\left(  k^{2}+\Delta\right)  \right)  \mathrm{d}g\\
&  =\left(  \left(  1-\chi\right)  k^{\mathbf{2}}+\Delta\right)  C+\left(
1-\chi\right)  \left(  k^{2}+\Delta\right)  \mathrm{d}g
\end{align*}
since%
\[
\left(  \Delta+k^{2}\right)  G=\delta_{0}%
\]
where $\delta_{0}$ is the delta-function on $\mathbb{R}^{3}$. Dividing by
$1-\chi$ we obtain
\[
\theta\mathbf{E}=\left(  k^{2}+\left(  1+\theta\right)  \Delta\right)
C+\left(  k^{2}+\Delta\right)  \mathrm{d}g,
\]
where $\theta\doteqdot\chi/\left(  1-\chi\right)  \ $hence%
\[
\theta\left(  \mathbf{E}-\Delta C\right)  =\left(  k^{2}+\Delta\right)
C+\left(  k^{2}+\Delta\right)  \mathrm{d}g.
\]
Differential operators $\mathrm{d}$ and $k^{2}+\Delta$ commute, therefore
\[
\mathrm{d}\left(  \theta\left(  \mathbf{E}-\Delta C\right)  \right)
\equiv\mathrm{d}\theta\wedge\left(  \mathbf{E}-\Delta C\right)  +\theta
\mathrm{d}\left(  \mathbf{E}-\Delta C\right)  =\left(  k^{2}+\Delta\right)
\mathrm{d}C\mathbf{.}%
\]
Multiplying by form $\mathbf{E}-\Delta C$ we kill the term with $\mathrm{d}%
\theta:$%
\[
\left(  \mathbf{E}-\Delta C\right)  \wedge\theta\mathrm{d}\left(
\mathbf{E}-\Delta C\right)  =\left(  \mathbf{E}-\Delta C\right)  \wedge\left(
k^{2}+\Delta\right)  \mathrm{d}C.
\]
By (\ref{13})$\ \mathrm{d}C=\mathrm{d}A_{0}=B^{\ast},\ $hence%
\begin{align*}
\Delta C  &  =\mathrm{d}^{\ast}\mathrm{d}C=\mathrm{d}^{\ast}B^{\ast}%
=\ast\mathrm{d}B,\\
\mathrm{d}\Delta C  &  =\mathrm{dd}^{\ast}\mathrm{d}C=\mathrm{dd}^{\ast
}B^{\ast}=\Delta B^{\ast}%
\end{align*}
\ which\ yields%

\[
\theta=\frac{\left(  \mathbf{E}-\Delta C\right)  \wedge\left(  k^{2}%
+\Delta\right)  \mathrm{d}C}{\left(  \mathbf{E}-\Delta C\right)
\wedge\mathrm{d}\left(  \mathbf{E}-\Delta C\right)  }=\frac{\left(
\mathbf{E}-\ast\mathrm{d}B\right)  \wedge\left(  k^{2}+\Delta\right)  B^{\ast
}}{\left(  \mathbf{E}-\ast\mathrm{d}B\right)  \wedge\left(  \mathrm{d}%
\mathbf{E}-\Delta B^{\ast}\right)  }.
\]
Finally
\[
\chi=\frac{1}{\theta+1}-1=-\frac{\left(  \mathbf{E}-\ast\mathrm{d}B\right)
\wedge\left(  k^{2}+\Delta\right)  B^{\ast}}{\left(  \mathbf{E}-\ast
\mathrm{d}B\right)  \wedge\left(  \mathrm{d}\mathbf{E}+k^{2}B^{\ast}\right)  }%
\]
and (\ref{9}) follows.$\ \blacktriangleright$

\section{Acquisition geometries}

Any acquisition geometry should be rich enough to guarantee\ reconstruction of
field $B$ on $\Omega.$ Let $x,y,z$ be an euclidean positively oriented
coordinate system on the physical space such that the static magnetic field
has the form $B_{0}=\left(  0,0,b\right)  ,$ $b>0$ in this system. The field
\
\begin{equation}
\mathit{B}^{+}=\left(  B_{x}+iB_{y}\right)  /2 \label{0}%
\end{equation}
is called \textit{positively rotating} part of $B\ $or \textit{transmit}
field. Several methods (sequences) are known that provide positively rotating
part of a magnetic field for ex. \cite{4}. According to \cite{3} the
"negatively rotating" part $\mathit{B}^{-}=\left(  B_{x}-iB_{y}\right)  /2$
can not be determined in this way.\textbf{ }The asymmetric role of coordinates
$x$ and $y$ in (\ref{0}) is defined with respect to $z$ coordinate which means
that the coordinate system $x,y,z$ is positively oriented. Note that the
orientation is indispensable feature of the Maxwell system which is
illustrated by Maxwell's right hand rule. It follows that such measurements of
the transmit field can be made for any positively oriented coordinate system
that is in the system obtained by rotation of the system $x,y,z.$\ 

\textbf{1. }The simple acquisition geometry for magnetic field $B$ is to fix
the body on the bed and rotate both around the central axis $y$ that is
parallel to the bed$\mathrm{.}$ Let $x,y,z$ be the laboratory euclidean
positively oriented (left-handed) system of coordinates,$\ $and let $\xi
,\eta,\zeta$ be the positively oriented system of coordinates that are
constant on the bed and the body when rotating that is%
\[
\xi=\cos\varphi x+\sin\varphi z,\eta=y,\ \zeta=-\sin\varphi x+\cos\varphi z
\]
where $\varphi$ is the rotation angle see Fig.1.%

\begin{figure}[ptb]%
\centering
\includegraphics[
natheight=7.499600in,
natwidth=13.333700in,
height=5.1742in,
width=9.1783in
]%
{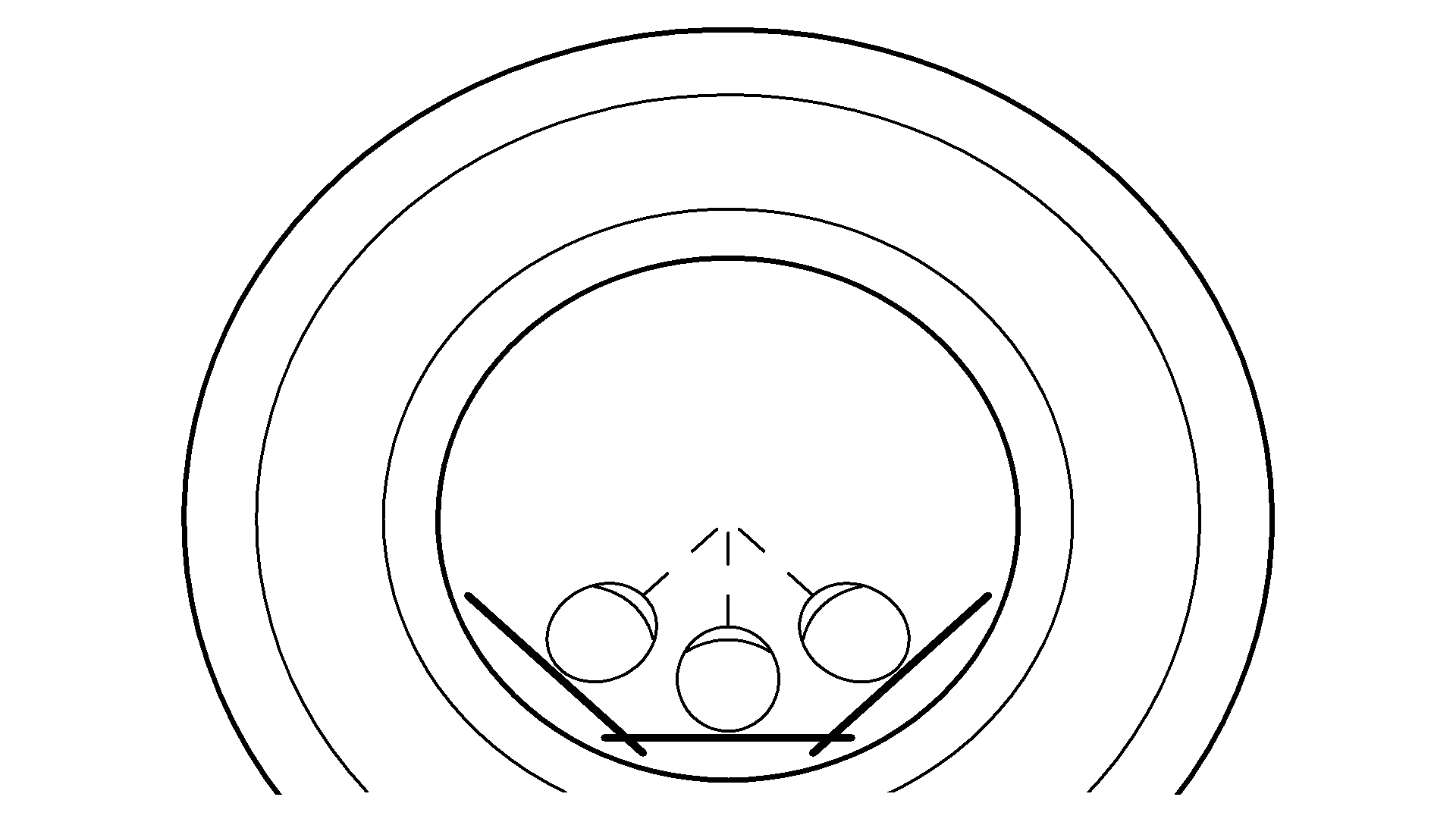}%
\end{figure}

Write magnetic field $B\left(  \varphi\right)  =\beta\ $by means of both
coordinate systems.
\begin{align}
B\left(  \varphi\right)   &  =B_{x}\left(  \varphi\right)  \mathrm{d}%
x+B_{y}\left(  \varphi\right)  \mathrm{d}y+B_{z}\left(  \varphi\right)
\mathrm{d}z,\label{5}\\
\beta &  =\beta_{\xi}\mathrm{d}\xi+\beta_{\eta}\mathrm{d}\eta+\beta_{\zeta
}\mathrm{d}\zeta\nonumber
\end{align}
and note that the body components of $\beta_{\xi},\beta_{\eta},\beta_{\zeta}$
of this field do not depend on $\varphi.\ $Keeping in view that
\[
\mathrm{d}\xi=\cos\varphi\mathrm{d}x+\sin\varphi\mathrm{d}z,\ \mathrm{d}%
\eta=\mathrm{d}y,\ \mathrm{d}\zeta=-\sin\varphi\mathrm{d}x+\cos\varphi
\mathrm{d}z
\]
we obtain%
\begin{align*}
B_{x}\left(  \varphi\right)   &  =\cos\varphi\beta_{\xi}-\sin\varphi
\beta_{\zeta},\ B_{y}\left(  \varphi\right)  =\beta_{\eta},\\
2\mathit{B}^{+}\left(  \varphi\right)   &  =B_{x}\left(  \varphi\right)
+iB_{y}\left(  \varphi\right)  =\cos\varphi\beta_{\xi}-\sin\varphi\beta
_{\zeta}+i\beta_{\eta}.
\end{align*}
Suppose that measurements of the transmit field$\ $are available for
$\varphi=-\psi,0,\psi\ $and consider system of linear equations%

\begin{align}
\cos\psi\beta_{\xi}-\sin\psi\beta_{\zeta}+i\beta_{\eta}  &  =2\mathit{B}%
^{+}\left(  \psi\right)  ,\nonumber\\
\cos\psi\beta_{\xi}+\sin\psi\beta_{\zeta}+i\beta_{\eta}  &  =2\mathit{B}%
^{+}\left(  -\psi\right)  ,\label{6}\\
\beta_{\xi}+i\beta_{\eta}  &  =2\mathit{B}^{+}\left(  0\right)  .\nonumber
\end{align}
The solution exists and\ is unique for any $\psi\neq0,\pi:$%
\begin{align*}
\beta_{\zeta}  &  =\frac{1}{\sin\psi}\left(  \mathit{B}^{+}\left(
-\psi\right)  -\mathit{B}^{+}\left(  \psi\right)  \right)  ,\ \beta_{\xi
}=\frac{1}{\cos\psi-1}\left(  \mathit{B}^{+}\left(  \psi\right)
+\mathit{B}^{+}\left(  -\psi\right)  -\mathit{B}^{+}\left(  0\right)  \right)
,\\
\beta_{\eta}  &  =\frac{i}{\cos\psi-1}\left(  \mathit{B}^{+}\left(
\psi\right)  +\mathit{B}^{+}\left(  -\psi\right)  +\left(  1-2\cos\psi\right)
\mathit{B}^{+}\left(  0\right)  \right)
\end{align*}
hence magnetic form $\beta$ is uniquely reconstructed by (\ref{5}). Now we can
apply (\ref{9}) to the fields $B=B\left(  0\right)  =\beta$ and $\mathbf{E.}$

\textbf{2. }Another geometry for scanning a body is as follows: the bed with
the body rotating in a tilted plane around the center $O$ in $B_{0}$
field.\ Let $\alpha,\ 0<\alpha<\pi/2$ be the constant tilting angle and
$e_{3}=\left(  0,-\sin\alpha,\cos\alpha\right)  .\ $For any $\varphi
,\ 0\leq\varphi<2\pi,$the vectors
\[
e_{1}=\left(  \cos\varphi,\sin\varphi\cos\alpha,\sin\varphi\sin\alpha\right)
,\ e_{2}=\left(  -\sin\varphi,\cos\varphi\cos\alpha,\cos\varphi\sin
\alpha\right)
\]
and $e_{3}\ $form the orthogonal frame since $e_{1}$ and $e_{2}$ belong to the
plane\ $P\ $through the origin orthogonal to $e_{3}$.\ This frame is a
positively oriented.\ Functions%
\[
\xi=\left\langle \left(  x,y,z\right)  ,e_{1}\right\rangle ,\ \eta
=\left\langle \left(  x,y,z\right)  ,e_{2}\right\rangle ,\ \zeta=\left\langle
\left(  x,y,z\right)  ,e_{3}\right\rangle
\]
are euclidean coordinates that are constant on $P$. To express the transmit
field (\ref{0}) by means of the body coordinates we write%
\begin{align*}
\mathrm{d}\xi &  =\cos\varphi\mathrm{d}x+\sin\varphi\cos\alpha\mathrm{d}%
y+\sin\varphi\sin\alpha\mathrm{d}z,\ \\
\mathrm{d}\eta &  =-\sin\varphi\mathrm{d}x+\cos\varphi\cos\alpha
\mathrm{d}y+\cos\varphi\sin\alpha\mathrm{d}z,\ \\
\mathrm{d}\zeta &  =-\sin\alpha\mathrm{d}y+\cos\alpha\mathrm{d}z
\end{align*}
and substitute in (\ref{5}). This gives%
\[
B_{x}\left(  \varphi\right)  =\cos\varphi\beta_{\xi}-\sin\varphi\beta_{\eta
},\ B_{y}\left(  \varphi\right)  =\sin\varphi\cos\alpha\beta_{\xi}+\cos
\varphi\cos\alpha\beta_{\eta}-\sin\alpha\beta_{\zeta}%
\]
and%
\[
2\mathit{B}^{+}\left(  \varphi\right)  =\cos\varphi\beta_{\xi}-\sin
\varphi\beta_{\eta}+i\left(  \sin\varphi\cos\alpha\beta_{\xi}+\cos\varphi
\cos\alpha\beta_{\eta}-\sin\alpha\beta_{\zeta}\right)  .
\]
Evaluating the transmit field for $\varphi=\varphi_{1},\varphi_{2},\varphi
_{3}$ we obtain the system like (\ref{6})\ for unknown functions $\beta_{\xi
},\beta_{\eta},\beta_{\zeta}$ that do not depend on $\varphi.$The determinant
of this\ system equals%

\[
\det=4\sin^{3}\alpha\sin\varphi_{12}\sin\varphi_{23}\sin\varphi_{31},\
\]
where $\varphi_{ij}=\left(  \varphi_{i}-\varphi_{j}\right)  /2.$\ It does
vanish if $\alpha\neq0,\pi$ and $\varphi_{ij}\neq0.\ $It follows that for
arbitrary $\alpha\neq0,\pi$, arbitrary different angles $\varphi_{i},\ $
$,\ $functions $\beta_{\xi},\beta_{\eta},\beta_{\zeta}$ are uniquely
determined\textrm{ }form data of fields$\ \mathit{B}^{+}\left(  \varphi
_{i}\right)  .$

\textbf{Algorithm for computation of Contrast source function}

1. Make three measurements of the transmit field according to geometry
\textbf{1} or \textbf{2}.

2. Calculate$\ $the field $\beta$ in the coordinate system of the body. and
set $B=B\left(  0\right)  =\beta.$

3. Apply (\ref{9}) to the fields $\beta$ and $\mathbf{E}$ and calculate the
quotient as the function $\chi\ $of the body coordinates$\ \xi,\eta,\zeta.\ $

\end{document}